\documentclass[12pt]{article}

\usepackage{amssymb,amsmath,exscale}

\usepackage{graphicx}

\usepackage{color}
\definecolor{marin}{rgb}   {0.,   0.3,   0.7} 
\definecolor{rouge}{rgb}   {0.8,   0.,   0.} 
\definecolor{sepia}{rgb}   {0.8,   0.5,   0.} 

\newtheorem{lemma}{Lemma}[section]

\newtheorem{theorem}[lemma]{Theorem}
\newtheorem{proposition}[lemma]{Proposition}

\newtheorem{remark}[lemma]{Remark}
\newtheorem{example}[lemma]{Example}	
\newtheorem{hypothesis}[lemma]{Hypothesis}
\newtheorem{notation}[lemma]{Notation}
\newtheorem{definition}[lemma]{Definition}
\newtheorem{conclusion}[lemma]{Conclusion}
\numberwithin{equation}{section}
\newcommand{\QED}{\mbox{}\hfill \raisebox{-0.2pt}{\rule{5.6pt}{6pt}\rule{0pt}{0pt}} 
          \medskip\par}             
\newenvironment{Proof}{\noindent
    \parindent=0pt\abovedisplayskip = 0.5\abovedisplayskip
    \belowdisplayskip=\abovedisplayskip{\bfseries Proof. }}{\QED}
\newenvironment{Proofof}[1]{\noindent
    \parindent=0pt\abovedisplayskip = 0.5\abovedisplayskip
    \belowdisplayskip=\abovedisplayskip{\bfseries Proof of #1. }}{\QED}

\newcommand{\dd}{\mathrm{d}}

\newcommand{\jb}{{\boldsymbol{j}}}

\newcommand{\kb}{{\boldsymbol{k}}}

\newcommand{\N}{\mathbb{N}}

\newcommand{\Mc}{\mathcal{M}}

\newcommand{\Kc}{\mathcal{K}}

\newcommand{\R}{\mathbb{R}}

\newcommand{\C}{\mathbb{C}}

\newcommand{\T}{\mathbb{T}}
\newcommand{\Z}{\mathbb{Z}}

\newcommand{\Zc}{\mathcal{Z}}

\newcommand{\Norm}[2]{\|#1\|\left.\vphantom{T_{j_0}^0}\!\!\right._{#2}}         
\newcommand{\SNorm}[2]{|#1|\left.\vphantom{T_{j_0}^0}\!\!\right._{#2}}             

\title{Sparse spectral approximations for computing polynomial functionals}        
\author{Erwan Faou, Fabio Nobile and Christophe Vuillot}      
\begin{document}
\maketitle
\abstract{
We give a new fast method for evaluating sprectral approximations of nonlinear polynomial functionals. We prove that the new algorithm is convergent if the functions considered are smooth enough, under a general assumption on the spectral eigenfunctions that turns out to be satisfied in many cases, including the Fourier and Hermite basis. 
\\[2ex]
}


\section{Introduction}

The goal of this paper is to introduce and analyze  a new method to compute spectral approximations of polynomial functionals  typically arising in spectral numerical methods applied to nonlinear partial differential equations. 

To describe the method and results, let us consider for instance the functional 
\begin{equation}
\label{eq:1}
X(u)(x) :=  u(x)^p
\end{equation}
where $u(x)$ is a smooth function on the one-dimensional torus $\T$ and $p \geq 2$ an integer. We can expand $u(x)$ as the Fourier series
$$
u(x) = \sum_{k \in \Z} u_k e^{ikx},
$$
where the $u_k \in \C$ are the Fourier coefficients associated with $u$. 
In this case, the functional  $X(u)(x) = \sum_{k \in \Z}e^{ikx} X_k(u)$ satisfies the convolution formula
\begin{equation}
\label{eq:2}
\forall\, k \in \Z, \quad X_k(u) = \sum_{\substack{(j_1,\ldots,j_p) \in \Z^p  \\ k = j_1 + \cdots + j_p }} u_{j_1} \cdots u_{j_p}. 
\end{equation}
To compute a numerical approximation of  such a quantity, a direct method would be prohibitive: if $u$ is approximated by $N$ coefficients, the sum on the right-hand side involves $N^{p-1}$ terms making the computational cost prohibitive for large $N$. 
That is why standard methods use the Fast Fourier Transform (FFT) algorithm to evaluate \eqref{eq:1} on grid points and an inverse FFT to go back to approximated Fourier coefficients $X_k$.  Though this method has the disadvantage to introduce aliasing problems due to the structure of FFT, it is very cheap in the sense that 
if the grid is made of $N$ points (and $u$ approximated by $N$ frequencies) the computational cost is of order $N \log N$. 

In many other situations like Hermite spectral methods, the problem is much harder because of the lack of fast transformation algorithm from collocation grid points to spectral variables (see however \cite{Arieh} for recent results by A. Iserles on a fast algorithm to compute Legendre coefficients). 

In this paper, we would like to show how a direct {\em sparse} approximation of \eqref{eq:2} of the form 
\begin{equation}
\label{eq:3}
\forall\, k \in \Z, \quad X_k^N(u) = \sum_{\substack{k = j_1 + \cdots + j_p \\ |j_1| \cdots |j_p| \leq N}} u_{j_1} \cdots u_{j_p}. 
\end{equation}
yields a consistent approximation of $X_k$ in the sense that we can control the difference $\Norm{X - X^N}{}$ in some Banach algebra, provided the function $u$ is sufficiently smooth. 

The big advantage of the representation \eqref{eq:3} is that the sum on the right-hand side involves only $\mathcal{O}(N (\log N)^{p-1})$ terms making the direct approximation at the spectral level possible and efficient. 


To have an idea of why this method is valid, let us calculate directly the difference
$$
X_k(u) - X_k^N(u) = \sum_{\substack{k = j_1 + \cdots + j_p \\ |j_1| \cdots |j_p| > N}} u_{j_1} \cdots u_{j_p}.
$$
We can write 
$$
|X_k(u) - X_k^N(u)| \leq \frac{1}{N^s} \sum_{\substack{k = j_1 + \cdots + j_p \\ |j_1| \cdots |j_p| > N}} |j_1|^s |u_{j_1}| \cdots |j_p^s| |u_{j_p}|.
$$
and we immediatly obtain the bound
\begin{equation}
\label{eq:4}
\Norm{X(u) - X^N(u)}{\ell^1} := \sum_{k \in \Z}|X_k(u) - X_k^N(u)|  \leq \frac{1}{N^s} \Big( \sum_{k \in \Z} |j|^s |u_j| \Big)^p = \frac{1}{N^s} \Norm{u}{\ell_s^1}^p. 
\end{equation}
Note that here we used $\ell^1$-based spaces (Wiener algebras) as they are the simplest to deal with polynomial nonlinearities in spectral representation. Even if similar results can be obtained for standard $\ell^2$-based spaces, we will state our result in these Banach spaces to avoid too many technical details. As high regularity $\ell^2$ and $\ell^1$ spaces are imbricated, this does not affect the validity of our approximation results. 

We see however that the previous proof does not extend straightforwardly to the case of Hermite functions. In this case, if $u(x)$ is defined on the real line $\R$ and decomposes into 
$$
u(x) = \sum_{k \in \N} u_k \chi_k(x)
$$
where the $\chi_j(x)$, $j \geq 0$ are normalized Hermite functions, then the Hermite coefficients $X_k$ of the product \eqref{eq:1} are given by  
\begin{equation}
\label{eq:pour}
\forall\,k \in \N,\quad 
X_k(u) = \sum_{(j_1,\ldots,j_p) \in \N^p  } a_{k;j_1,\ldots,j_p}u_{j_1} \cdots u_{j_p},
\end{equation}
where 
\begin{equation}
\label{eq:prodchi}
a_{k;j_1,\ldots,j_p} = \int_{\R} \chi_k(x) \chi_{j_1}(x) \cdots \chi_{j_p}(x) \dd x
\end{equation}
are the integrals of products of Hermite functions. Note that in this situation, the coefficients are non zero even in the case where $k \neq j_1 + \cdots + j_p$. To obtain a convergence result similar to \eqref{eq:4} we thus see that we need a non trivial control of these coefficients. To this aim, we take advantage of the recent work by B. Gr\'ebert, R. Imekraz and E. Paturel, see \cite{GIP}, in which bounds are given for the coefficients $a_{k;j_1,\ldots,j_p}$ that allow to prove that $X$ acts on high-regularity Sobolev spaces.  Note that this Hermite case is of particular importance because of the lack of fast Hermite transform, while in practice Hermite spectral methods are quite natural and widely used in many applications fields like Bose-Einstein condensate simulations and Fokker-Planck equations. 

In Section 2, we give a very general result in an abstract setting by assuming explicit bounds on the coefficients $a_{k;j_1,\ldots,j_p}$ in \eqref{eq:pour}. This result covers the case of Fourier and Hermite basis, spherical harmonics functions, and eigenfunctions of operators of the form $-\Delta + V$ in dimension one with Dirichlet or periodic boundary conditions. 

To cover different situations, we introduce general sparse sets of indices of the form $|k|^\alpha |j_1| \cdots |j_p| \leq N$ where $\alpha = 0$ or $1$. 

In the case where the {\em momentum} $k - j_1 -\cdots j_p$ is bounded in the sum defining  $X_k$ (like in the Fourier case, see \eqref{eq:3}), the set of non zero coefficients will be indeed of size $\mathcal{O}(N (\log N)^{p-1})$ for $\alpha = 0$. However in more general situations like Hermite approximation, the set in $k$ and $(j_1,\ldots,j_d)$ will be of size $\mathcal{O}(N^{2} (\log N)^{p-1})$ for $\alpha = 0$ (if $|k| < N$) and $\mathcal{O}(N (\log N)^p)$ for $\alpha = 1$. The effect of this parameter $\alpha$ is only a slight deterioration of the rate of convergence of the approximation, but it reduces drastically the computational cost of the method for large $N$ in the Hermite case.  

In Section 3, we show how an iterative implementation of the algorithm yields a convergent approximation of the product of $p$ functions with a cost of order $\mathcal{O}(p N \log N)$ instead of $\mathcal{O}(N (\log N)^{p-1})$. We give an error estimate for this case as well. 

In Section 4 we detail the case of periodic exponential functions (the Fourier basis) and discuss the possible extensions to eigenfunctions of operators of the form $-\Delta + V$. In Section 5 we consider the Hermite case and show by numerical experiments that the error bounds are optimal. 

\section{An abstract result}
We consider $\Zc = \Z^d$ or $\N^d$ for $d \geq 1$.
For $u = (u_j)_{j \in \Zc} \in \C^{\Zc}$ we set 
\begin{equation}
\label{eq:defls}
\Norm{u}{\ell_s^1} = \sum_{j \in \Zc} \Norm{j}{}^s |u_j|, 
\end{equation}
where $\Norm{j}{} = \max(1,|j^1|, \ldots, |j^d|)$ for $j = (j^1,\ldots, j^d) \in \Zc$. 
We also define the norm 
\begin{equation}
\label{eq:defl2}
\Norm{u}{\ell_s^2} = \Big(\sum_{j \in \Zc} \Norm{j}{}^{2s} |u_j|^2\Big)^{\frac12},
\end{equation}
and using the Cauchy-Schwartz inequality, we can easily prove that if $s' - s > d/2$, there exists a constant $C$ such that for all $u$, we have  
\begin{equation}
\label{eq:imbric}
\Norm{u}{\ell_s^2} \leq  \Norm{u}{\ell^1_s} \leq C \Norm{u}{\ell_{s'}^2}. 
\end{equation}

For a given integer $p \geq 2$, 
we aim at approximating a function $X: (\ell^1_s)^p \to \ell_s^1$ defined by $X(u^1,\ldots,u^p) = (X_\ell(u^1,\ldots,u^p))_{\ell \in \Zc}$ where 
\begin{equation}
\label{eq:Xu}
\forall\, \ell \in \Zc,\quad 
X_\ell(u^1,\ldots,u^p) = \sum_{j_1,\cdots j_p \in \Zc^p} a_{\ell ; j_1 \cdots j_p} u^1_{j_1} \cdots u^p_{j_p}, 
\end{equation}
with given coefficients $a_{\ell ; j_1 \cdots j_p} \in \C$. 
We use the following notation: for a multiindex $\jb = (j_1,\ldots,j_p)$ and $\ell\in \Zc$, we define the {\em momentum} 
\begin{equation}
\label{eq:defmom}
\Mc(\ell,\jb) = \ell - j_1 - \cdots - j_p. 
\end{equation}
We will also sometime use the notation $a_{\ell;\jb}$ to denote the coefficient $a_{\ell;j_1 \cdots j_p}$. 

\subsection{Sparse sets of frequencies}

We consider a subset $\Kc \subset \Zc$. We will typically consider the case where $\Kc = \Zc$, a bounded set of $\Zc$ or a sparse set of indices of $\Zc$. We assume that $\Kc$ is equipped with a function $|\cdot |$ measuring the size of multi-indices of the form  $j = (j^1,\ldots,j^d) \in \Zc$. We assume that there exist positive constants $c_0$, $C_0$ and $\sigma$ such that 
\begin{equation}
\label{eq:leq}
\forall\, j \in \Zc, \quad c_0\Norm{j}{} \leq |j| \leq C_0 \Norm{j}{}^\sigma. 
\end{equation}
We then set for $u \in \C^\Zc$ (compare \eqref{eq:defls})
\begin{equation}
\label{eq:deflsS}
\SNorm{u}{\ell_s^1} = \sum_{j \in \Zc} |j|^s |u_j|, 
\end{equation}
and using \eqref{eq:leq} we obtain
\begin{equation}
\label{eq:imbricspar}
c\Norm{u}{\ell_s^1} \leq \SNorm{u}{\ell_s^1}  \leq  C\Norm{u}{\ell_{\sigma s}^1}  
\end{equation}
for some constant $c$ and $C$ independent of $u$.
As particular cases of application, we mainly have in mind the two following situations: 
\begin{itemize}
\item[(i)] $\Kc$ is a set of the form 
\begin{equation}
\label{eq:passpard}
\Kc_M = \{ \, j \in \Zc \, | \, |j| \leq M\, \}, \quad \mbox{with}\quad |j| := \Norm{j}{}, 
\end{equation}
where $M \in \overline \N$ can be equal to $+\infty$ in which case $\Kc_M = \Zc$. In this situation, we have $C_0 = c_0 = \sigma = 1$ in the inequality \eqref{eq:leq}. Note that for a given $M$, we have $\sharp \Kc_M \leq C M^d$ for some constant $C$ independent on $M$, where $\sharp F$ denotes the cardinal of the set $F$. 
\item[(ii)] $\Kc$ is a sparse set of the form 
\begin{equation}
\label{eq:spard}
\Kc_M^* = \{ \,j \in \Zc \, | \, |j| \leq M\,\}, \quad \mbox{where}\quad |j| := \prod_{n = 1}^{d} (1 + |j^n|), 
\end{equation}
for some given $M \in \overline\N$. In this case, using the inequality of arithmetic and geometric means, \eqref{eq:leq} is valid with $\sigma = d$. In this situation, we have $\sharp \Kc_M^* \leq C M (\log M)^{d-1}$ for some constant $C$ independent of $M$ (see for instance \cite{Griebel,Zenger}). 
\end{itemize}

For a fixed $\alpha \in \{0,1\}$ and $N \geq 0$, we define the following approximation $X^{N,\alpha}(u) = (X^{N,\alpha}_\ell)_{\ell \in \Kc}$ of $X(u)$: 
\begin{equation}
\label{eq:approx}
\forall\, \ell \in \Kc,\quad 
X_\ell^{N,\alpha}(u^1,\ldots,u^p) = \sum_{\substack{j_1,\cdots j_p \in \Kc^p \\ |\ell|^\alpha |j_1| \cdots |j_p| \leq N}} a_{\ell ; j_1 \cdots j_p} u^1_{j_1} \cdots u^p_{j_p}. 
\end{equation}

The next Lemma estimates the number of non zero terms involved in the definition of $X^{N,\alpha}(u)$ in the two cases (i) and (ii) described above. 

%

\begin{lemma}
\label{lem:size}
The cardinals of the sparse sets of indices can be estimated as follows: Let $\alpha \in \{0,1\}$ and $p \geq 1$. There exists a constant $C$ depending only on $d$ and $p$ such that, for all $M$ and $N \geq 1$, we have 
\begin{itemize}
\item[(i)] With $\Kc_M$  defined by \eqref{eq:passpard}, and $|j| = \Norm{j}{}$ for all $j \in \Zc$, then  
$$
\sharp \{ \ell, j_1,\cdots j_p \in \Kc_M^p \, | \,  |\ell|^\alpha|j_1| \cdots |j_p| \leq N \}  \leq C (\sharp \Kc_M)^{(1 - \alpha)}N^d (\log N)^{p - 1 + \alpha }, 
$$
with the convention $(\sharp \Kc_M)^{0} = 1$ when $\Kc_M = \Zc$, that is $M = +\infty$.
\item[(ii)] With $\Kc_M^*$ and $|j|$ the sparse norm defined by \eqref{eq:spard}, then we have 
\begin{multline*}
\sharp \{ \ell, j_1,\cdots j_p \in (\Kc_M^*)^p \, | \,  |\ell|^\alpha |j_1| \cdots |j_p| \leq N \}  \\
\leq C (\sharp \Kc_M)^{(1 - \alpha)}  N (\log N)^{d(p+\alpha)- 1}.  
\end{multline*}
\end{itemize}

\end{lemma}

\begin{Proof} 
The proof of (ii) is classical (see for instance \cite{Griebel,Zenger}) using the fact that in this case, $|j| = \prod_{k = 1}^d ( 1 + |j^k|)$ when $j = (j^1,\ldots,j^d)$, so that 
$$
|\ell|^\alpha |j_1| \cdots |j_p| = \Big( \prod_{k = 1}^d ( 1 + |\ell^k|) \Big)^{\alpha} \prod_{k = 1}^d \prod_{n = 1}^p ( 1 + |j^k_n|). 
$$
which yields the result for $\alpha = 1$ (independently on $M$). The case $\alpha = 0$ is treated similarly. 

The proof of (i) is a consequence of the fact that  for all $N \geq 1$ and $p \geq 1$, 
$$
 \sharp \{ j_1,\cdots ,j_p \in \Zc^p \, | \,  \Norm{j_1}{} \cdots \Norm{j_p}{} \leq N \} \leq C_p N^d (\log N)^{p-1}. 
$$
for some constant $C_p$ depending on $p$ and $d$. We prove this by induction on $p$: for $p = 1$ the result is clear using $\Norm{j}{} = \max(1,|j^1|, \ldots, |j^d|) \in \N\backslash\{0\}$ for $j = (j^1,\ldots,j^d) \in \Zc$. Let us assume that it holds for $p - 1 \geq 1$. We have 
\begin{equation*}
\begin{split}
\sharp \{ j_1, & \cdots ,j_p \in  \Zc^p \, | \,  \Norm{j_1}{} \cdots \Norm{j_p}{} \leq N \} 
\\
&= \sum_{ k = 1}^N \sharp \{ j_1,\cdots j_{p-1} \in \Zc^{p-1} \, | \,  \Norm{j_1}{} \cdots \Norm{j_{p-1}}{} \leq \frac{N}{k} \}  \times \sharp \{ j \in \Zc\, | \, \Norm{j}{} = k\, \}, \\
&\leq 2^d \sum_{ k = 1}^N C_{p-1} \big(\frac{N}{k}\big)^d  \big(\log \frac{N}{k}\big)^{p-2} \times d k^{d-1} \\
&\leq 2^d d C_{p-1} N^d (\log N)^{p-2} \sum_{ k = 1}^N \frac{1}{k} \leq C_p N^d (\log N)^{p-1} 
\end{split}
\end{equation*}
for some constant $C_p$ depending on $p$ and $d$. This yields the result. Here we used the fact that we calculate explicitly that for $k \geq 2$, $\sharp \{ j \in \Zc\, | \, \Norm{j}{} = k\, \} = d k^{d-1}$, while for $k = 1$, this number is equal to $2^d$, with the definition of $\Norm{j}{}$. 
\end{Proof}

As we will see below, the previous result can be refined when the coefficients $a_{\ell,\jb}$ in \eqref{eq:Xu} have some special structure implying a decay property with respect to the momentum $\Mc(\ell,\jb)$ defined in \eqref{eq:defmom}. We will consider theses cases more in detail in the section devoted to the Fourier case. 

\subsection{Error estimate}

The goal of this section is to give an estimate of the error 
$$
\Norm{X(u^1,\ldots,u^p) - X^{N,\alpha}(u^1,\ldots,u^p)}{\ell_s^1} 
$$
for smooth $u^1, \ldots, u^p$ and where $X^{N,\alpha}$ is defined by \eqref{eq:approx} for some given $\alpha\in \{0,1\}$ and $N\geq 1$.


We make a general hypothesis on the coefficients $a_{\ell;\jb}$ involved in the definition of the functional $X$.  

\begin{definition}
Let $\kb = (k_1,\ldots, k_q) \in \Zc^q$  with $q \geq 1$ a multi-index. For $n = 1,\ldots,q$, we set $\mu_n(\kb)$ the $n$-th largest integer amongst $\Norm{k_1}{},\ldots,\Norm{k_q}{}$, so that we have $\mu_1(\kb) \geq \mu_2(\kb) \geq \mu_3(\kb) \geq \cdots$. 
\end{definition}

%
We make the following hypothesis: 
\begin{hypothesis}
\label{hyp:1}
 There exist $\nu\geq 0$, $\theta \in [0,1]$ such that for all $R$, there exists $c_R$ such that for all $\ell \in \Zc$, and all $\jb  = (j_1,\ldots, j_p) \in \Zc^p$,  we have 
\begin{equation}
\label{eq:bound}
|a_{\ell;\jb}| \leq c_R \mu_3(\kb)^\nu \Big(
\frac{\mu_2(\kb)^\theta \mu_3(\kb)^{1 - \theta}}{\mu_2(\kb)^\theta \mu_3(\kb)^{1 - \theta} + \mu_1(\kb) - \mu_2(\kb)} \Big)^R.
\end{equation}
where  $\kb =  (\ell, \jb) = (\ell,j_1,\ldots,j_p)$. 
\end{hypothesis}

Let us make some comments on this definition. Such bounds (with $\theta = 0$) were used in several recent works \cite{Delort1,Delort2,Zoll,Bam03,BG06,Greb07} to prove long time existence results on nonlinear PDEs set on manifolds with different kind of boundary conditions (compact manifold,  Dirichlet, etc...). It holds true in many situations where the $a_{\ell,\jb}$ are products of the form \eqref{eq:prodchi} with functions $\chi_k$ defining a $L^2$ Hilbert basis on a manifod $M$, like the Fourier basis on a torus. It is also valid (with $\theta = 0$) in the case of spherical harmonics, see \cite{Delort1,Delort2}, and when $\chi_k$ are {\em well localized with respect to the exponentials}, see \cite{Bam03,BG06} and Definition 5.3 of \cite{Greb07}. This last situation corresponds to the case where the $\chi_k$ are eigenfunctions of an operator $-\Delta + V$ with Dirichlet boundary conditions in dimension 1, and with a smooth periodic potential $V$. 

More recently this was extended to Hermite functions basis diagonalizing the quantum harmonic oscillator operator, see \cite{GIP}. In this case the previous bound holds true but for $\theta = 1/2$.

The main result of this section is the following. 
\begin{theorem}
\label{th1}
Assume that the coefficients $a_{\ell;\jb}$ of the function $X(u^1,\ldots,u^p)$  satisfy the Hypothesis \ref{hyp:1} for some constants $\nu \geq 0$ and $\theta \in [0,1]$, and let $X^{N,\alpha}$ be the approximation \eqref{eq:approx} defined for $\alpha \in \{0,1\}$, $N \geq 1$ and $(\Kc,|\cdot|) \subset (\Zc,\Norm{\cdot}{})$ satisfying \eqref{eq:leq} for some constant $\sigma \geq 1$. Let  $\kappa > d$ be fixed. Then for all $s\geq 0$ and $s' \geq \max(\sigma s + \theta \kappa, (1-\theta)\kappa + \nu)$, 
there exists a constant $C$ such that for all $N$ and for all functions $u^i \in \ell^1_{\sigma s'}$, $i = 1,\ldots, p$,  we have the estimate
\begin{equation}
\label{eq:th}
\Norm{X(u^1,\ldots,u^p) - X^{N,\alpha}(u^1,\ldots,u^p)}{\ell_s^1} \leq C N^{-\beta(s,s')} \prod_{i = 1}^{p} \SNorm{u^i }{\ell_{s'}^1}, 
\end{equation}
where 
\begin{equation}
\label{eq:beta}
\beta(s,s')  = \min\Big(  \frac{s' -  \sigma s -  \theta\kappa}{\sigma\alpha + 1}, s' - (1 - \theta) \kappa - \nu \Big).  
\end{equation}
\end{theorem}

To prove this Theorem, we will use the following technical Lemma. 
The proof of this Lemma is postponed to the Appendix.

\begin{lemma}
\label{lem:1}
Let  $\kappa > d$. 
Assume that $a_{\ell;\jb}$ satisfies the previous Hypothesis \ref{hyp:1} for some constants $\nu \geq 0$ and $\theta \in [0,1]$. Then for all $r \geq 0$, there exists a constant $C_r$ such that for all $\jb =  (j_1,\ldots,j_p) \in \Zc^p$, 
\begin{equation}
\label{eq:lem1}
\sum_{\ell \in \Zc} \Norm{\ell}{}^r |a_{\ell;\jb}| \leq C_r \mu_1(\jb)^{r + \theta \kappa} \mu_2(\jb)^{(1 - \theta) \kappa + \nu}. 
\end{equation}

\end{lemma}

\begin{Proofof}{Theorem \ref{th1}}
We set for $\ell\in \Zc$, 
\begin{eqnarray*}
R_\ell(u^1,\ldots,u^p) &=& X_\ell(u^1,\ldots,u^p) - X_\ell^{N,\alpha}(u^1,\ldots,u^p)\\
&=&\sum_{\substack{j_1,\cdots j_p \in \Zc^p \\  |\ell|^\alpha|j_1| \cdots |j_p| > N}} a_{\ell ; j_1 \cdots j_p} u^1_{j_1} \cdots u^p_{j_p} . 
\end{eqnarray*}
For some $t \leq s'$, we can write 
\begin{eqnarray*}
\SNorm{R}{\ell_s^1} &\leq& \sum_{\substack{\ell,j_1,\cdots j_p \in \Zc^{p+1} \\  |\ell|^\alpha|j_1| \cdots |j_p| > N}} |\ell|^{s} | a_{\ell ; j_1 \cdots j_p} u^1_{j_1} \cdots u^p_{j_p}| \\
&\leq& \frac{1}{N^{s' - t }}\sum_{\substack{\ell,j_1,\cdots j_p \in \Zc^p \\  |\ell|^{\alpha }|j_1| \cdots |j_p| > N}}  \frac{|\ell|^{s + \alpha ( s' - t)} |a_{\ell ; j_1 \cdots j_p}|}{  |j_1|^t \cdots |j_p|^t}\,  |j_1|^{s'} |u^1_{j_1}| \cdots |j_p|^{s'}|u^p_{j_p}|. 
\end{eqnarray*}
Hence we get using \eqref{eq:leq}
$$
\Norm{R}{\ell_s^1}  \leq \frac{1}{c_0} \SNorm{R}{\ell_s^1} \leq \frac{C(t) }{N^{s'-t}} \prod_{i = 1}^p \SNorm{u^i}{\ell_{s'}^1},
$$
where 
$$
C(t) := \frac{C_0}{c_0} \sup_{(j_1,\ldots, j_p) \in \Zc^p} \frac{1}{  |j_1|^t \cdots |j_p|^t}\sum_{\ell} \Norm{\ell}{}^{\sigma s + \sigma \alpha ( s' - t)} |a_{\ell ; j_1 \cdots j_p}|, 
$$
where $C_0$ is the constant appearing in \eqref{eq:leq}. 
Applying the previous Lemma with $r = \sigma s + \sigma \alpha (s' - t)$ and using again \eqref{eq:leq} we see that $C(t)$ will be finite if 
$$
t = \max ( \sigma s + \sigma\alpha ( s' - t) + \theta \kappa , (1 - \theta) \kappa + \nu). 
$$
or equivalently
$$
t = \max(\frac{\sigma s + \sigma \alpha s' + \theta \kappa}{\sigma\alpha + 1}, (1 - \theta) \kappa + \nu), 
$$
in which case $s' - t = \beta(s,s')$. This shows the result. 
\end{Proofof}

\section{Iterative approximations}\label{sec:iterative}

We consider now the case where $X(u^1,\ldots,u^p)$ corresponds to the product operator of $p$ functions $u^i = \sum_{j \in \Zc} u^i_j \chi_j(x)$, $1 \leq i \leq p$, where $\chi_j(x)$ is an orthonormal basis of $L^2(M)$ where $M$ is a manifold (typically $M = \T^d$ or $\R^d$). In this case, the coefficients $a_{\ell;j_1,\ldots,j_n}$ are given by the integrals
$$
a_{\ell;j_1,\ldots,j_n} = \int_{M} \chi_\ell \chi_{j_1} \cdots \chi_{j_n} \dd M. 
$$
We will see in the example below that bound \eqref{eq:bound} holds in many situations such as the Fourier basis on $\T^d$ and the Hermite basis on $\R^d$. In such a case, we identify a function $u$ with its coefficients $u_j$ and talk about $u \in \ell_s^1$ by a slight abuse of notation. 

In the previous section, we have proven that for two functions $u^1$ and $u^2$, the function $X^{N,\alpha}(u^1,u^2)$ yields a good approximation of the product $u^1 u^2 = X(u^1,u^2)$ if these functions are smooth. Now for three functions $u^1$, $u^2$ and $u^3$, instead of approximating the product $u^1u^2 u^3$ by using $X^{N,\alpha}(u^1,u^2,u^3)$, which generates a computational cost of order $\mathcal{O}(N (\log N)^{3})$ in dimension $d = 1$ and for $\alpha = 1$ (see Lemma \ref{lem:size}), we might use the following algorithm: 
\begin{enumerate}
\item Compute the approximation $v = X^{N,\alpha}(u^1,u^2)$ of the product $u^1u^2$
\item Compute $X^{N,\alpha}(v,u^3)$ as approximation of $u^1u^2u^3$. 
\end{enumerate}
In other words, we replace $X^{N,\alpha}(u^1,u^2,u^3)$ by $X^{N,\alpha}(X^{N,\alpha}(u^1,u^2),u^3)$.

Obviously  the cost of this algorithm is of order $\mathcal{O}(2 N (\log N)^2)$ for $\alpha = 1$, instead of $\mathcal{O}(N (\log N)^{3})$ (in dimension $d = 1$, see Lemma \ref{lem:size}). Such an iterative approximation can be easily generalized to any product of $p$ functions, and the global cost is of order $\mathcal{O}(p N (\log N)^2)$ for $\alpha = 1$, instead of $\mathcal{O}(N (\log N)^{p})$. As we will see now, an error estimate of the same kind as in the previous section remains valid for such sparse approximations. For simplicity, we only present the result in the case where $|\cdot| = \Norm{\cdot}{}$, which implies  $\Norm{u}{\ell_s^1} = \SNorm{u}{\ell_s^1}$. 

This is given by the following result: 

\begin{theorem}
\label{th22}
Let $u^i(x)$, $i \geq 0$ be given functions. For all $i \geq 0$, let us define the functions $U^{N,\alpha}(u^1,\ldots,u^i)$ by induction as follows: $U^{N,\alpha}(u^1) = u^1$ and for $i \geq 1$, 
$$
U^{N,\alpha}(u^1,\ldots,u^{i+1}) =  X^{N,\alpha}(U^{N,\alpha}(u^1,\ldots,u^{i}), u^{i+1}),
$$
where $\alpha \in \{0,1\}$ is fixed and $X^{N,\alpha}$ defined in \eqref{eq:approx} for $(\Kc,|\cdot|) \subset (\Zc,\Norm{\cdot}{})$ where $|j| = \Norm{j}{}$ for all $j \in \Zc$. 
Then for all $p$, the function $U^{N,\alpha}(u^1,\ldots,u^{p})$ is an approximation of $X(u^1, \ldots u^p) = u^1 \cdots u^p$ in the following sense: Assume that the coefficients $a_{\ell;\jb}$ satisfy the Hypothesis \ref{hyp:1} for some constants $\nu \geq 0$ and $\theta \in [0,1]$, and let  $\kappa > d$ be fixed. Then for all $p \in \N$,  $s\geq 0$ and $s' \geq \max(s + (p-1) \theta \kappa, (1 - \theta)\kappa + \nu)$, 
there exists a constant $C$ such that for all $N$  we have the estimate
\begin{equation}
\label{eq:th2}
\Norm{X(u^1,\ldots,u^p) - U^{N,\alpha}(u^1,\ldots,u^p)}{\ell_s^1} \leq C N^{-\beta_p(s,s')} \prod_{i = 1}^{p} \Norm{u^i }{\ell_{s'}^1}, 
\end{equation}
where 
\begin{equation}
\label{eq:beta2}
\beta_p(s,s')  = \min\Big(  \frac{s' -  s - (p-1)\theta\kappa}{\alpha + 1}, s' - (p-3)\theta \kappa -  \kappa - \nu )\Big).  
\end{equation}

\end{theorem}
\begin{Proof}
As $N$ and $\alpha$ are fixed, we set $U^{i} := U^{N,\alpha}(u^1,\ldots,u^i)$. 
For $p = 2$, the estimate is the one given in Theorem \ref{th1} with $\sigma = 1$. Assume that it holds for $p - 1 \geq 2$. In particular, we have for all $s''\geq 0$ and $s'\geq \max(s'' + (p-1) \theta \kappa, (1 - \theta)\kappa + \nu)$
$$
\Norm{U^{p-1}}{\ell_{s''}^1} \leq \big( 1 +   C N^{-\beta_{p-1}(s'',s')} \big)\prod_{i = 1}^{p-1} \Norm{u^i }{\ell_{s'}^1}, 
$$
for some constant $C$ depending on $s'$, $s''$ and $p$. Here we use the fact that in the case where $|j| = \Norm{j}{}$ the norms $\Norm{\cdot}{\ell_s^1}$ and $\SNorm{\cdot}{\ell_s^1}$ coincide. 
Now using the definition of $U^p$, we can write 
\begin{multline}
\label{eq:hjk}
U^{p}  - X(u^1 ,\ldots, u^p) =  X^{N,\alpha}(U^{p-1}, u^{p}) - U^{p-1} \cdot u^p \\
+ (U^{p-1} - X(u^1, \ldots, u^{p-1})) \cdot u^p . 
\end{multline}
As a direct consequence of Lemma \ref{lem:1}, we easily see that the following holds: for $s  > (1 - 2\theta) \kappa + \nu$, 
and for $u = \sum_{j \in \Zc} u_j \chi_j$ and $v = \sum_{j \in \Zc} v_j \chi_j$  in $\ell_s^1$, we have 
$$
\Norm{uv}{\ell_s^1} \leq \sum_{\ell,j_1,j_2 \in \Zc} \Norm{\ell}{}^s |a_{\ell;j_1 j_2}| |u_{j_1}| |v_{j_2}| \leq C_s \Norm{u}{\ell_{s+ \theta \kappa}^1} \Norm{v}{\ell_{s+\theta \kappa}^1}. 
$$
Using this inequality and \eqref{eq:hjk} we obtain for $s'' \geq \max(s +  \theta \kappa, (1 - \theta)\kappa + \nu)$, using \eqref{eq:th2} for $p = 2$, 
\begin{multline*}
\Norm{U^{p}  - X (u^1, \ldots, u^p)}{\ell_s^1}  \leq   C N^{- \beta_2(s,s'')} \Norm{U^{p-1}}{\ell_{s''}^1} \Norm{u^p}{\ell_{s''}  } \\
+ \Norm{U^{p-1} - X(u^1, \ldots, u^{p-1})}{\ell_{s+ \theta \kappa}^1}\Norm{u^p}{\ell_{s + \theta \kappa}^1}
\end{multline*}
and hence, for some constant $C$ depending on $s$, $s'$, $s''$ and $p$, 

\begin{multline*}
\Norm{U^{p}  - X(u^1, \cdots ,u^p)}{\ell_s^1}  \\
\leq  C \Big(  N^{- \beta_2(s,s'')} \big( 1 +    N^{-\beta_{p-1}(s'',s')} \big)+  N^{-\beta_{p-1}(s + \theta \kappa,s')} \Big)
 \times \prod_{i = 1}^p  \Norm{u^i}{\ell_{s'}^1}. 
\end{multline*}
We take $s'' = s' - (p-2)\theta \kappa$, so that $\beta_{p-1}(s'',s') = 0$. For this $s''$ we have 
\begin{eqnarray*}
\beta_2 (s,s'') &=& \min\Big(  \frac{s'' -  s - \theta\kappa}{\alpha + 1}, s'' +  \theta \kappa - \nu - \kappa \Big)\\
&=&\min\Big(  \frac{s' -  s - (p-1)\theta\kappa}{\alpha + 1}, s' - (p-3)\theta \kappa -  \kappa - \nu \Big)
\end{eqnarray*}
Moreover, we have 
$$
\beta_{p-1}(s + \theta \kappa,s') = \min\Big(  \frac{s' -  s - \theta \kappa - (p-2)\theta\kappa}{\alpha + 1}, s' - (p-4)\theta \kappa -  \kappa - \nu  \Big). 
$$
On taking the minimum between $\beta_{p-1}(s + \theta \kappa,s')$ and $\beta_2 (s,s'')$, we obtain the result. 
\end{Proof}

In the rest of this paper, we will show how this Theorem can be applied to many situations including the discretization of polynomials in Fourier or Hermite  basis. 

\section{Fourier basis}

We consider now functions $u(x)$ defined on $x \in \T^d$. We consider functionals of the form 
\begin{equation}
\label{eq:aup1}
X(u^1,\cdots,u^p)(x) = b(x) \, u^1(x) \cdots u^p(x),
\end{equation}
where $b(x)$ is a given function defined on the torus $\T^d$. With a function $u(x) \in \C$, $x = (x^1,\cdots x^d) \in \T^d$, and for a given $j = (j^1,\cdots,j^d) \in \Zc:= \Z^d$ we associate the Fourier coefficients 
$$
u_j = \frac{1}{(2\pi)^d} \int_{\T^d} u(x) e^{-i j \cdot x} \, \dd x,
$$
where $j \cdot x = j^1 x^1 + \cdots j^d x^d$. In this case, the coefficients $a_{\ell;j_1 \cdots j_p}$ defined in \eqref{eq:Xu} can be calculated explicitely, and for given $\ell \in \Zc = \Z^d$ and $\jb = (j_1,\ldots, j_p) \in \Zc^p$.   
\begin{equation}
\label{eq:fou}
a_{\ell;j_1 \cdots j_p} = 
\sum_{k \in \Z^d} b_{k} \frac{1}{(2\pi)^d}  \int_{\T^d} e^{i(- \ell + k + j_1 + \cdots + j_p)\cdot x} \dd x\\
=  b_{\Mc(\ell,\jb)}, 
\end{equation}
where the numbers $b_k$ are the Fourier coefficients associated with the function $b(x)$, and with the definition \eqref{eq:defmom} of the momentum $\Mc(\ell,\jb)$. Here we use the very special property of the exponential functions $e^{ij \cdot x}$ that the product of two basis functions is again a basis function. We assume that $b(x)$ extends to an  analytic function on a complex strip $U_\rho := \T^d \times i [-\rho,\rho]^d$ around the torus, which implies by standard Cauchy estimates that 
\begin{equation}
\label{eq:hypbk}
\forall\, k \in \Z^d, \quad |b_k| \leq  D e^{- \rho |k|}, 
\end{equation}
where $D = \sup_{z \in U_\rho} |b(z)| $. 

With this calculation, we can prove the following result: 
\begin{proposition}
If the Fourier coefficients $b_k$ of the function $b(x)$ satisfy the analytic estimate \eqref{eq:hypbk}, then the coefficients $a_{\ell,j_1\cdots j_p}$ defined in \eqref{eq:fou} satisfy the Hypothesis \ref{hyp:1} with $\nu = 0$ and $\theta = 0$. 


\end{proposition}

Hence we see that when $b$ is analytic, we will have $\beta(s,s') = \frac{s' - s}{\alpha + 1}$ in the formula \eqref{eq:beta}, provided $s'$ and $s$ are large enough. The proof of the previous proposition can be found in \cite{Bam03,Greb07}. As explained in these references, the same result holds true when the function $u(x)$ is decomposed on a Hilbert basis $e_j(x)$, $j \in \Zc^d$ that is {\em well-localized} with respect to the exponential. This includes in particular the case where $e_j(x)$ are the eigenfunctions of a differential operator of the form $u(x) \mapsto -\Delta u(x) + V(x) u(x)$ for some smooth periodic potential function $V(x)$ in dimension $d = 1$. We refer to \cite{Greb07} for extensive discussions on the subject. 

Let us mention that in the particular case where $b(x)$ is a trigonometric polynomial containing only a finite number of frequencies, the use of the parameter $\alpha = 1$ is not mandatory to obtain sparse set of indices. This is a consequence of the Lemma below:  

\begin{lemma}
\label{lem:size2}
Considering the approximation \eqref{eq:approx}, we assume that there exists $q \geq 0$ such that  
$$
|\Mc(\ell,\jb)| > q \Longrightarrow a_{\ell; \jb} = 0. 
$$
The cardinals of the sparse sets of indices with $\alpha = 0$ can be estimated as follows: Let $p \geq 1$, then there exists a constant $C$ depending only on $d$, $q$ and $p$ such that, for all $M$ and $N \geq 1$, we have 
\begin{itemize}
\item[(i)] With $\Kc_M$  defined by \eqref{eq:passpard}, and $|j| = \Norm{j}{}$ for all $j \in \Zc$, then  
$$
\sharp \{ \ell, j_1,\cdots j_p \in \Kc_M^p \, | \,  |j_1| \cdots |j_p| \leq N \}  \leq C   N^d (\log N)^{p - 1 }. 
$$
\item[(ii)] With $\Kc_M^*$ and $|j|$ the sparse norm defined by \eqref{eq:spard}, then we have 
$$
\sharp \{ \ell, j_1,\cdots j_p \in (\Kc_M^*)^p \, | \,   |j_1| \cdots |j_p| \leq N \}  
\leq C   N (\log N)^{dp- 1}.  
$$
\end{itemize}
\end{lemma}
\begin{Proof} 
For a fixed $\ell$ in the sets considered, the estimates can be obtained similarly as in the proof of Lemma \ref{lem:size}, and are independent of $\ell$. Now 
under the assumption on the momentum, if $a_{\ell; \jb} \neq 0$ then we have necessarily $\Mc(\ell,\jb) = m \in \Z^d$ with $|m| \leq q$. Summing in $m$ then yields the result (with a constant proportional to $q^d$). 
\end{Proof}

Note that the case considered in the introduction corresponds to $b(x) = 1$ and $q = 0$ in the previous Lemma.

\begin{figure}
\centering
\includegraphics[width=.45\textwidth]{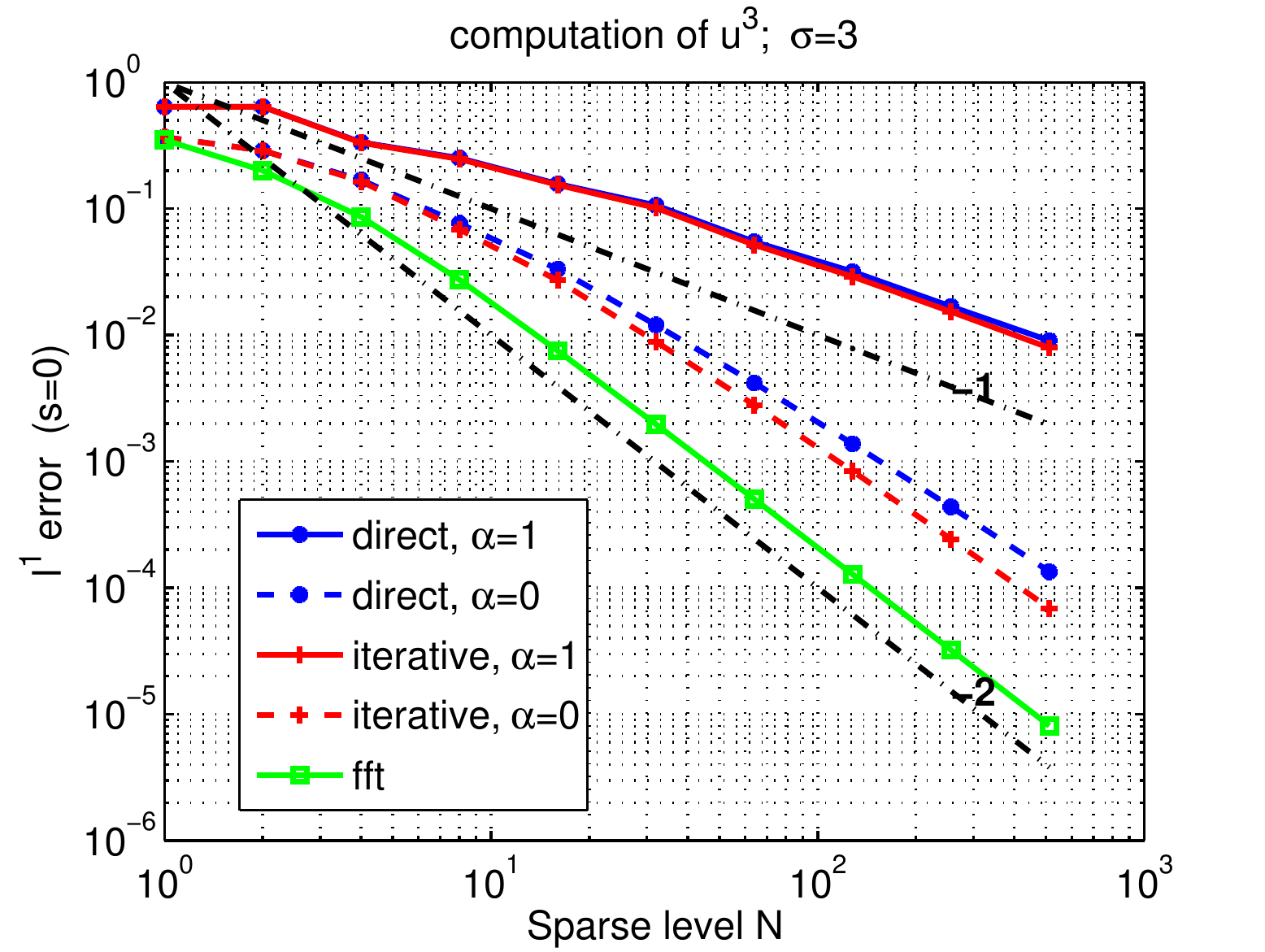} \quad
\includegraphics[width=.45\textwidth]{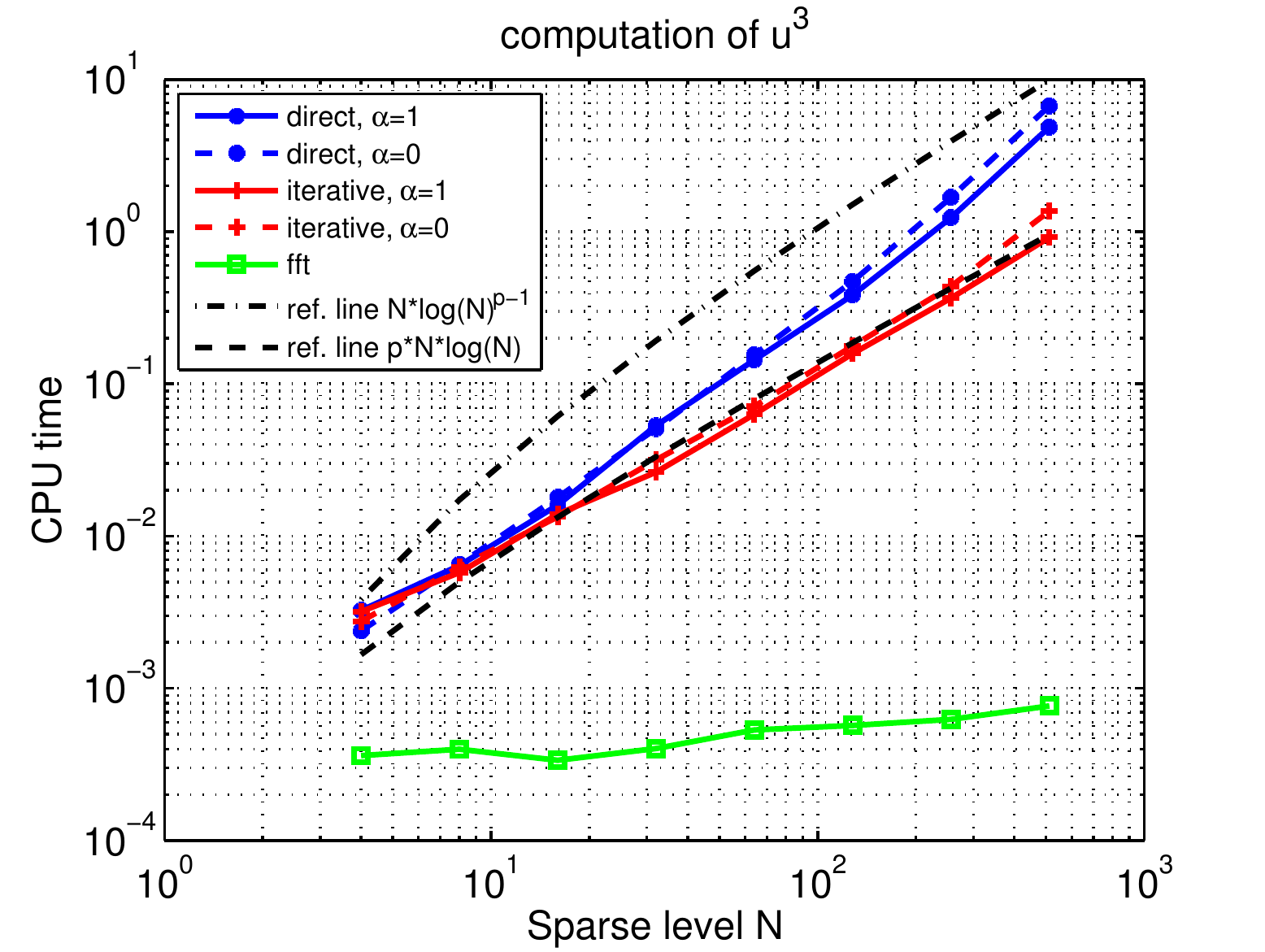}
\caption{Sparse approximation of $u^3$ for $\sigma=3$. Left:
  convergence of the $\ell^1$-error; Right: CPU time}
\label{fig:l1errorp3}
\end{figure}

We show now on a numerical example the accuracy of the estimates
above. We consider the function $u(x) = \sum_{k\in \Z} u_k e^{ikx}$
with $u_k=(1+|k|)^{-\sigma}$ so that $u \in \ell_{s'}^1$ for $s' <
\sigma-1$. We compute $u^p$ by the \emph{direct} method
\eqref{eq:approx} and the \emph{iterative} algorithm described in
Section \ref{sec:iterative}. In both cases we expect a maximal
convergence rate $O(N^{\frac{\sigma-1}{\alpha+1}})$ in $\ell^1$ (that
is for $s=0$). In figure \ref{fig:l1errorp3} (left) we plot in log
scale the $\ell^1$-error versus the sparse level $N$ in the case $p=3$
for the different approximation methods. 
In figure \ref{fig:l1errorp3} (right) we plot the estimated CPU time
together with the theoretical bounds $CN(\log N)^{p-1+\alpha}$ for the
direct method and $CpN(\log N)^{1+\alpha}$ for the iterative one. For
convenience we plot only the theoretical bounds for $\alpha=0$.  The
version $\alpha=0$ is clearly more accurate than $\alpha=1$. On the
other hand it has only a minimal extra cost, so for this particular
example it is clearly preferable. It should be pointed out, however,
that this is due to the very simple form of the functional $X(u)=u^3$
for which $a_{\ell; \jb} = 0$ if $|\Mc(\ell,\jb)|\ne 0$.

\begin{figure}
\centering
\includegraphics[width=.45\textwidth]{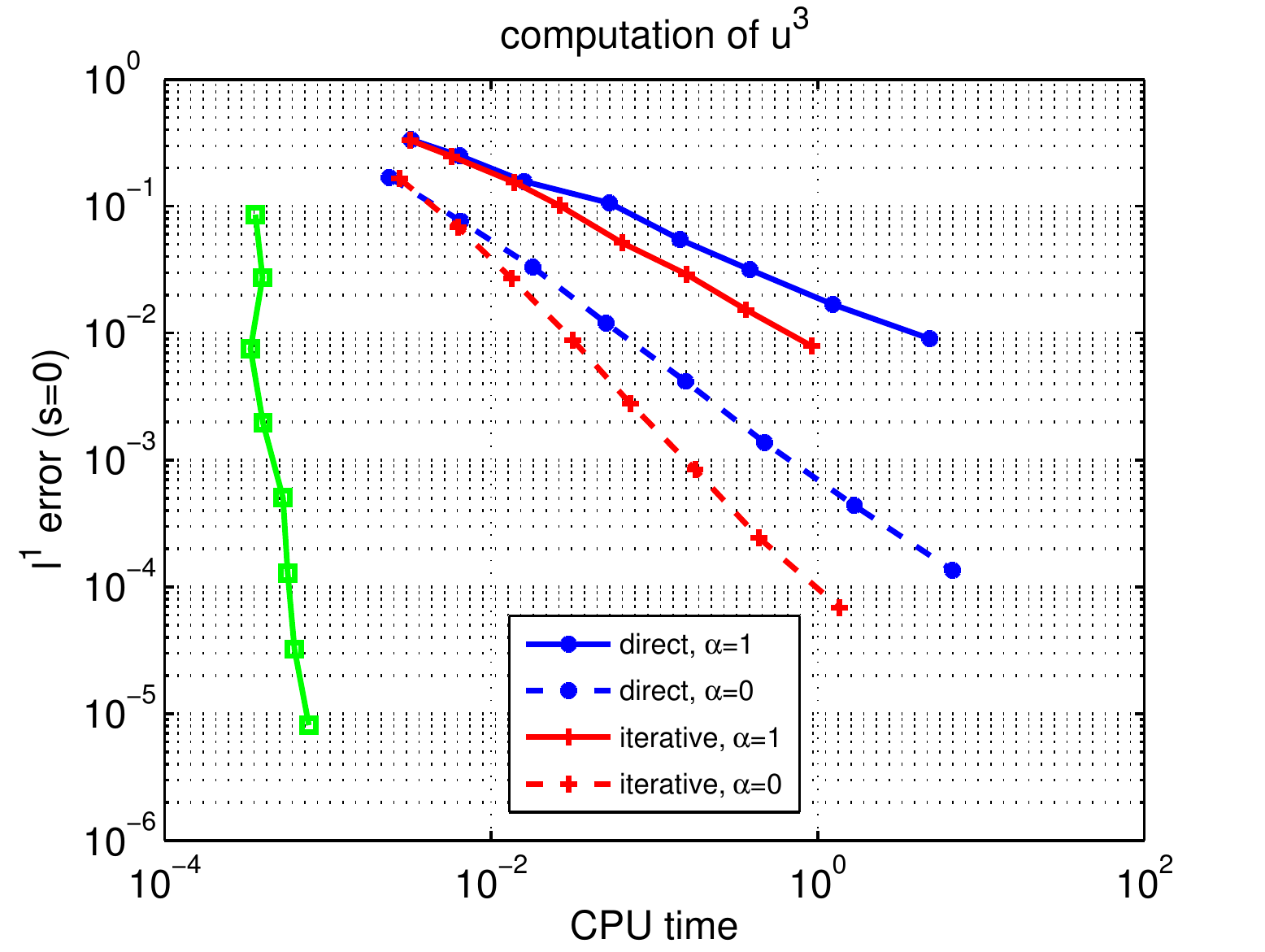}
\caption{Sparse approximation of $u^3$ for $\sigma=3$. $\ell^1$-error versus CPU time}
\label{fig:complexity}
\end{figure}

Figure \ref{fig:complexity} shows the error versus the CPU time. It is
clear from this plot the advantage of the iterative algorithm with
respect to the direct method, as well as the advantage of $\alpha=0$
with respect to $\alpha=1$.

\begin{figure}
\centering
\includegraphics[width=.45\textwidth]{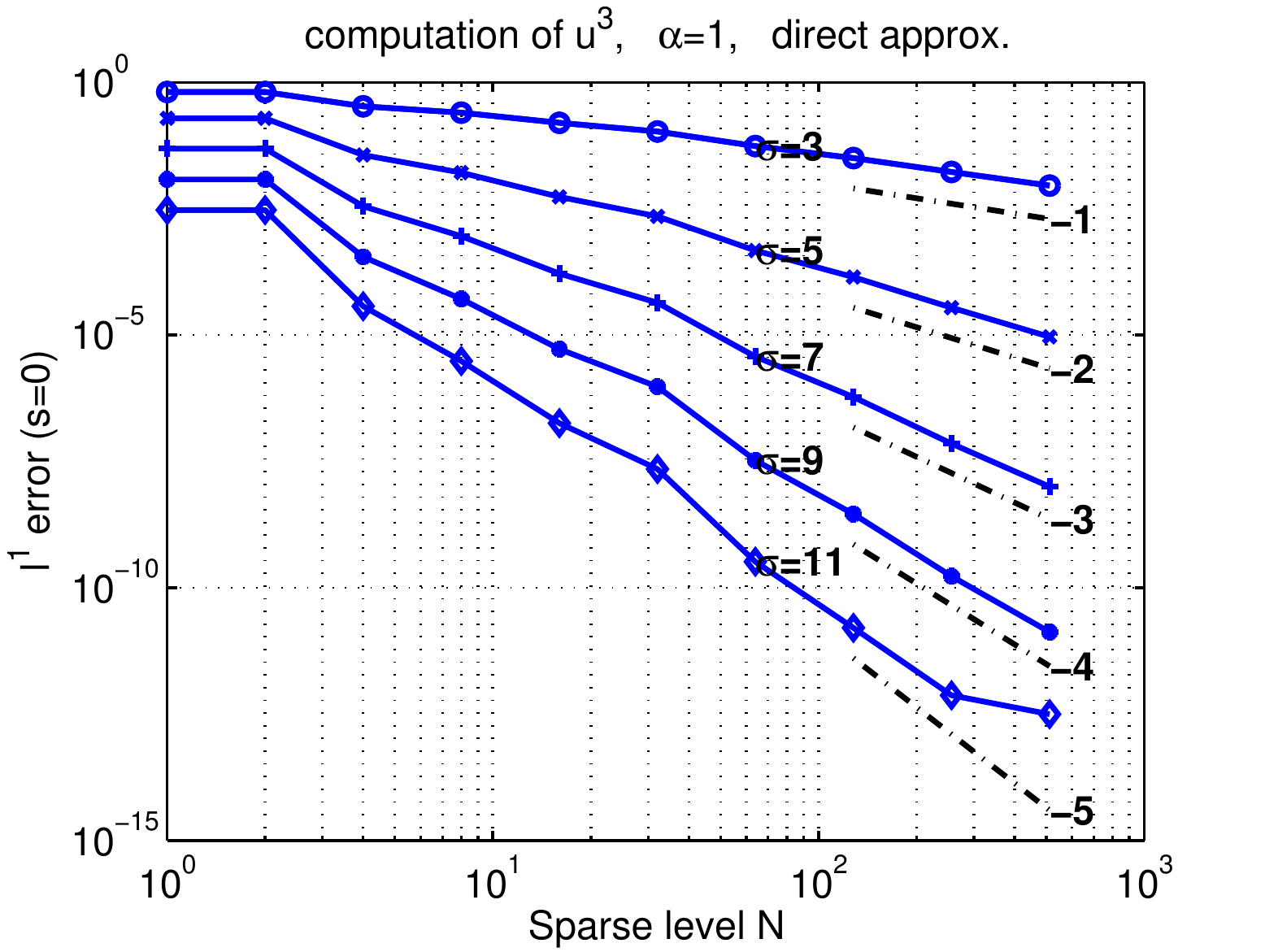}
\caption{Convergence of the sparse approximation of $u^3$ with the direct formula \eqref{eq:approx} and $\alpha=1$}
\label{fig:sigma}
\end{figure}

Finally, in figure \ref{fig:sigma} we show the convergence of the
$\ell^1$-error, still in the case $p=3$ but for different values of
$\sigma$. We consider here only the case of $\alpha=1$ and the direct
formula \eqref{eq:approx}. The results in the other cases are
analogous. For all values of $\sigma$ we recover the expected
theoretical rate of convergence.

\section{Hermite}

We consider now the case where $u(x)$ is defined on the real line ($x \in \R$) and the basis $(\chi_j)_{j \in \N}$ is given by 
the set of normalized Hermite functions defined by the formula
\begin{equation}
\label{eq:Tdef}
T \chi_j := - \frac{\dd^2\chi_j}{\dd x^2} (x) + x^2 \chi_j (x) = (2j  + 1) \chi_j(x), \quad j \in \N, 
\end{equation}
with the condition $\Norm{\chi_j}{L^2(\R)} = 1$. 
Note that here, with the notation of the previous sections, we have $\Zc = \N$. 
For all $j \in \N$, the Hermite functions are given by 
$$
\chi_n(x) = \frac{H_n(x) }{\sqrt{2^n n! \sqrt{\pi}}} e^{-x^2/2}
$$
where $H_n(x)$ is the $n$-th Hermite polynomial with respect with the weight $e^{-x^2}$. Recall that these Hermite polynomials satisfy 
$$
\forall\, n, m \in \N, \quad \int_{\R} H_n(x) H_m(x) e^{-x^2} \dd x = 2^n n! \sqrt{\pi} \delta_{nm},
$$
and the induction relations: 
$$
H_0(x) = 1, \quad\mbox{and}\quad  H_{n+1}(x) = 2x H_{n}(x) - 2n H_{n-1}(x), \quad n \geq 1. 
$$
In this situation, $(\chi_j(x))_{j \in \N}$ is a Hilbert basis of $L^2(\R)$ and for a given real function $u(x)$, we can write 
$$
u(x) = \sum_{j \in \N} u_j \chi_j(x),\quad \mbox{where} \quad u_j = \int_{\R} u(x) \chi_j(x)\, \dd x. 
$$

Here, note that the Hilbert space  associated with the norm $\ell_s^2$ defined in \eqref{eq:defl2} coincides with the domain of the operator $T^{s}$ (see \eqref{eq:Tdef}). Using standard notations, the classical space $\tilde H^s$  defined by 
$$
\tilde H^s = \{ u(x) \in H^s(\R)\, | \, x \mapsto x^p \partial_x^q u(x) \in L^2(\R)\, \quad \mbox{for}\quad 0 \leq p + q \leq s\,\}
$$
corresponds with the domain of the operator $T^{s/2}$ (see for instance \cite{Hel84}) and hence with $\ell_{s/2}^2$. In particular, we can write owing to \eqref{eq:imbric}
$$
c\Norm{u}{\tilde H^{s}} \leq \Norm{u}{\ell_{s/2}^1}\leq C \Norm{u}{\tilde H^{s'}},
$$
provided $s - s' > d$, and for some positive constants $c$ and $C$ independent of $u$.

Let us now consider the functional $X(u)(x) = u(x)^p$ for $p \in \N$. In this case, the coefficients $a_{\ell,j_1 \cdots j_p}$ in \eqref{eq:Xu} are given by the formula, for $(\ell,j_1,\ldots,j_p) \in \N^{p+1}$. 
\begin{equation}
\label{eq:aKL}
a_{\ell;j_1\ldots j_p} = \int_{\R} \chi_\ell(x) \chi_{j_1}(x) \cdots \chi_{j_p}(x) \dd x.  
\end{equation}
The following Proposition can be found in \cite{GIP}, Proposition 3.6: 
\begin{lemma}
For all $\nu > 1/8$ and all $R> 0$, there exists $c_R$ such that for all $\ell \in \N$, and all $\jb  = (j_1,\ldots, j_p) \in \N^p$,  we have 
\begin{equation}
\label{eq:boundherm}
|a_{\ell;\jb}| \leq c_R \frac{\mu_3(\kb)^\nu}{\mu_1(\kb)^{\frac{1}{24}}} \Big(
\frac{\mu_2(\kb)^{\frac12} \mu_3(\kb)^{\frac12}}{\mu_2(\kb)^{\frac12} \mu_3(\kb)^{\frac12} + \mu_1(\kb) - \mu_2(\kb)} \Big)^R.
\end{equation}
where  $\kb =  (\ell, \jb) = (\ell,j_1,\ldots,j_p)$. In particular, these coefficients satisfy
 Hypothesis \ref{hyp:1} with $\theta = 1/2$ and $\nu > 1/8$. 
\end{lemma}

Note that the estimate given in \cite{GIP} is slightly better than the one given in \ref{eq:bound} because of the presence of the term $\mu_1(\kb)^{\frac{1}{24}}$ in the denominator in \eqref{eq:boundherm}. 

\begin{remark}
In this paper, we will only consider the case of Hermite functions in dimension 1. The extension to higher dimension can be made using the framework of \cite{GIP}, Section 3.2.  
\end{remark}

In this situation, and when $\alpha = 1$, we obtain a convergence rate of order $N^{\frac{s - s' + \kappa}{2}}$ for $\kappa > 1$ and sufficiently large $s'$ for the algorithm described in Theorem \ref{th1}, and $N^{\frac{s - s' + (p-1)\kappa/2}{2}}$ for the iterative algorithm (see Theorem \ref{th22}). 

In the following we will illustrate these results by numerical simulation. In all the computations presented below, the coefficients \eqref{eq:aKL} are approximated in double machine precision by using Gauss-Hermite quadrature rules with the packages provided by J. Burkhardt\footnote{\href{http://people.sc.fsu.edu/~jburkardt/cpp\_src/hermite\_rule/hermite\_rule.html}{http://people.sc.fsu.edu/$\sim$jburkardt/cpp\_src/hermite\_rule/hermite\_rule.html}}.

We first consider the case where $p = 3$, and for given number $\sigma$, we consider the fonctions $u(x) = \sum_{n \geq 0} u_n \chi_n(x)$ with $u_n = (1 + n)^{-\sigma}$ so that $u \in \ell_{s'}^1$ for $s' < \sigma-1$. 
\begin{figure}[ht]
\begin{center}
\rotatebox{0}{\resizebox{!}{0.5\linewidth}{%
   \includegraphics{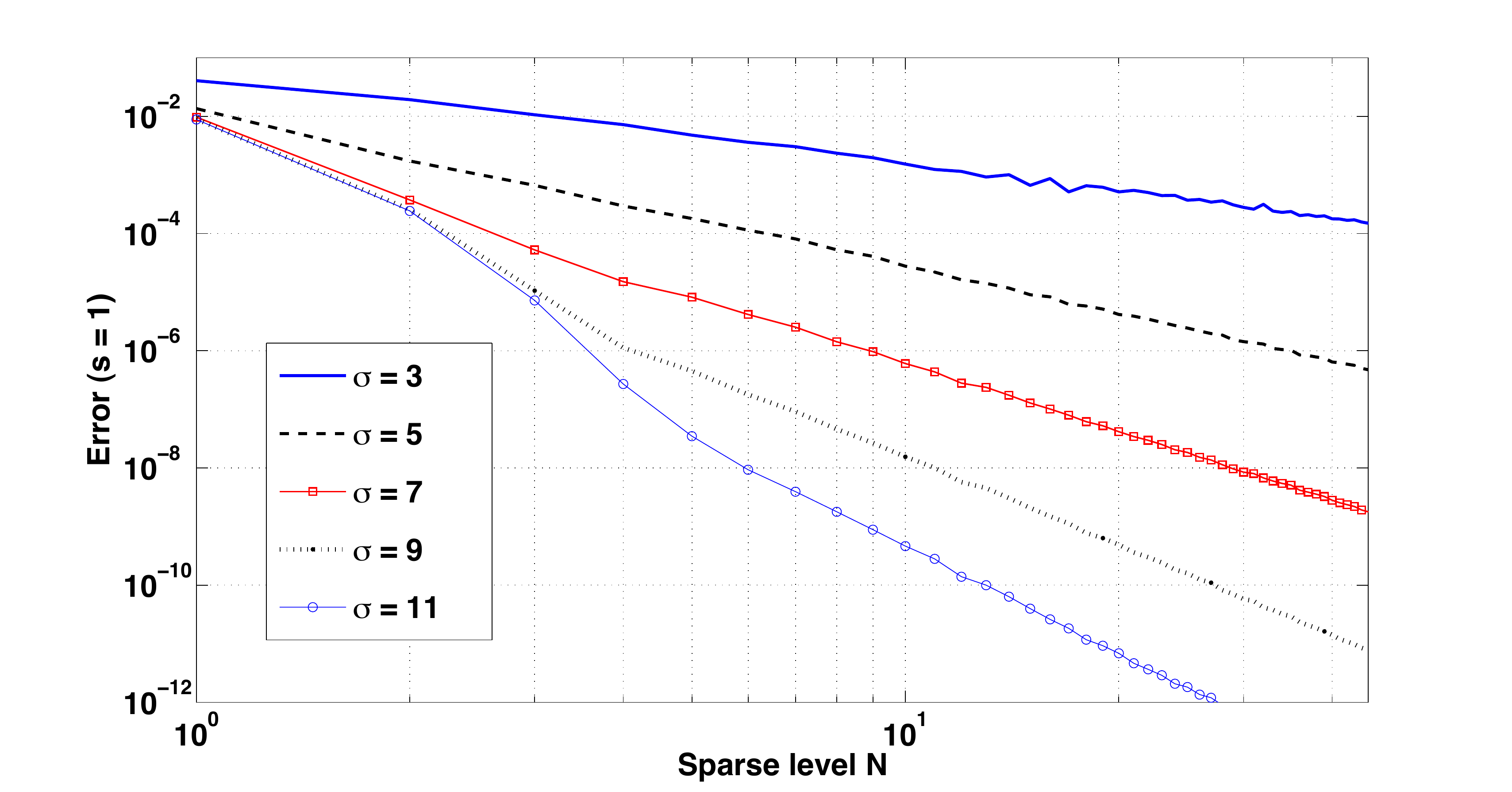}}}
\end{center}
\caption{Convergence of the sparse approximation}
\label{fig1}
\end{figure}
Hence in this case, we expect a maximal convergence rate of order $\mathcal{O}(N^{ - \frac{\sigma-1 - \kappa}{2}})$ in $\ell^1$ (that is for $s = 0$) and when $\alpha = 1$. In Figure \ref{fig1} we plot in log-log scale the error measured in  $\ell^1$ norm between the approximation $X^{N,1}(u,u,u)$ and the exact solution $u^3$ whose Hermite coefficients are approximated using a Hermite transform with 500 points. The convergence rates observed correspond to the theoretical estimate \eqref{eq:th}. 

\begin{figure}[ht]
\begin{center}
\rotatebox{0}{\resizebox{!}{0.5\linewidth}{%
   \includegraphics{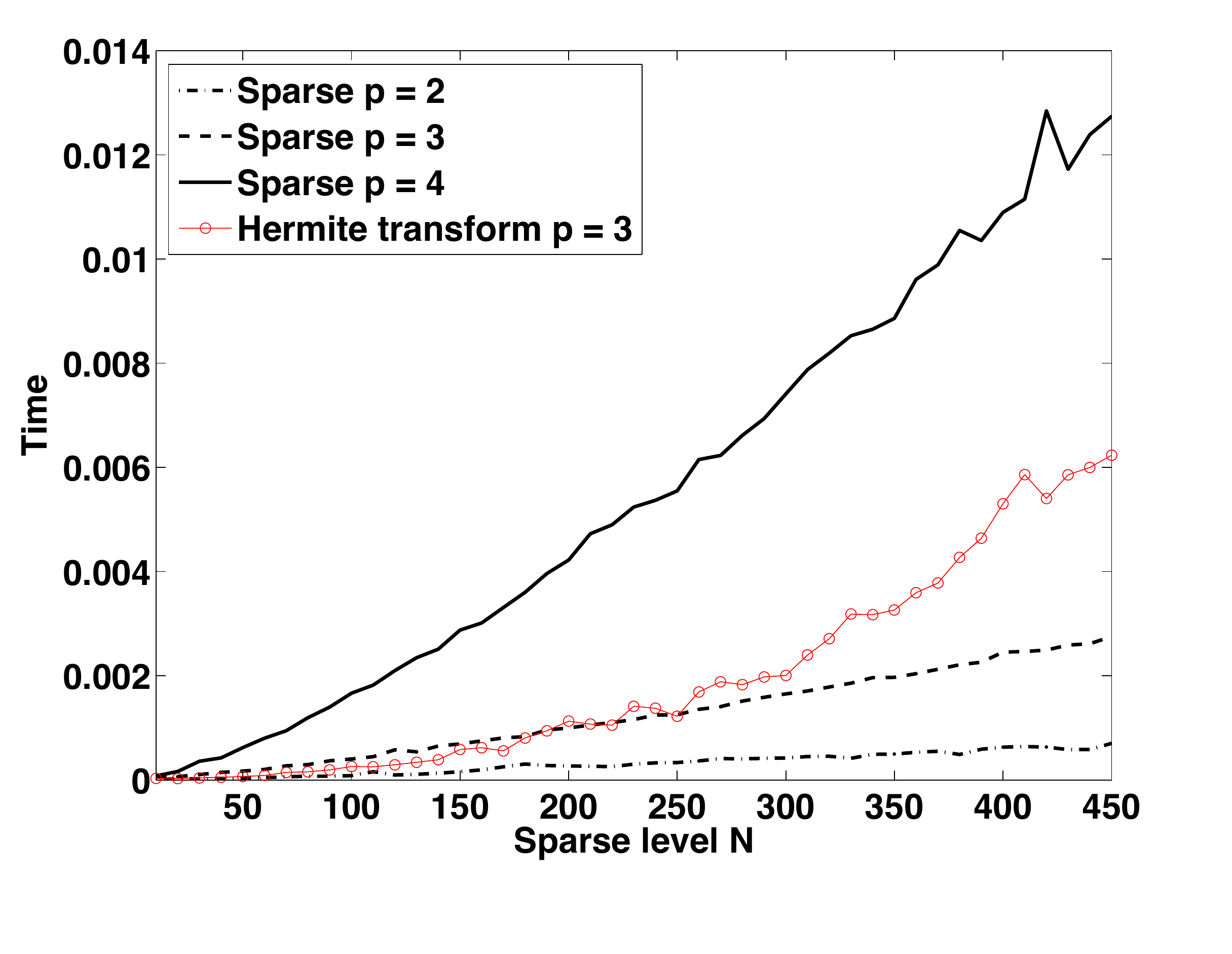}}}
\end{center}
\caption{Time versus sparse level $N$}
\label{fig2}
\end{figure}

In figure \ref{fig2}, we plot the time required by the algorithm in the cases $p = 2,3$ and $p= 4$, to compute the Hermite coefficients of $u(x)^p$. As expected, the time increases when $p$ becomes large, which is in accordance with the cost of order $N (\log N)^{p}$ predicted by Lemma \ref{lem:size}. We compare with the cost of the Hermite transform algorithm with $N$ points, which is in $\mathcal{O}(N^2)$. Note that for this latter method, the cost does not significantly differ with $p$, and only $p = 3$ is shown. 

\begin{figure}[ht]
\begin{center}
\rotatebox{0}{\resizebox{!}{0.5\linewidth}{%
   \includegraphics{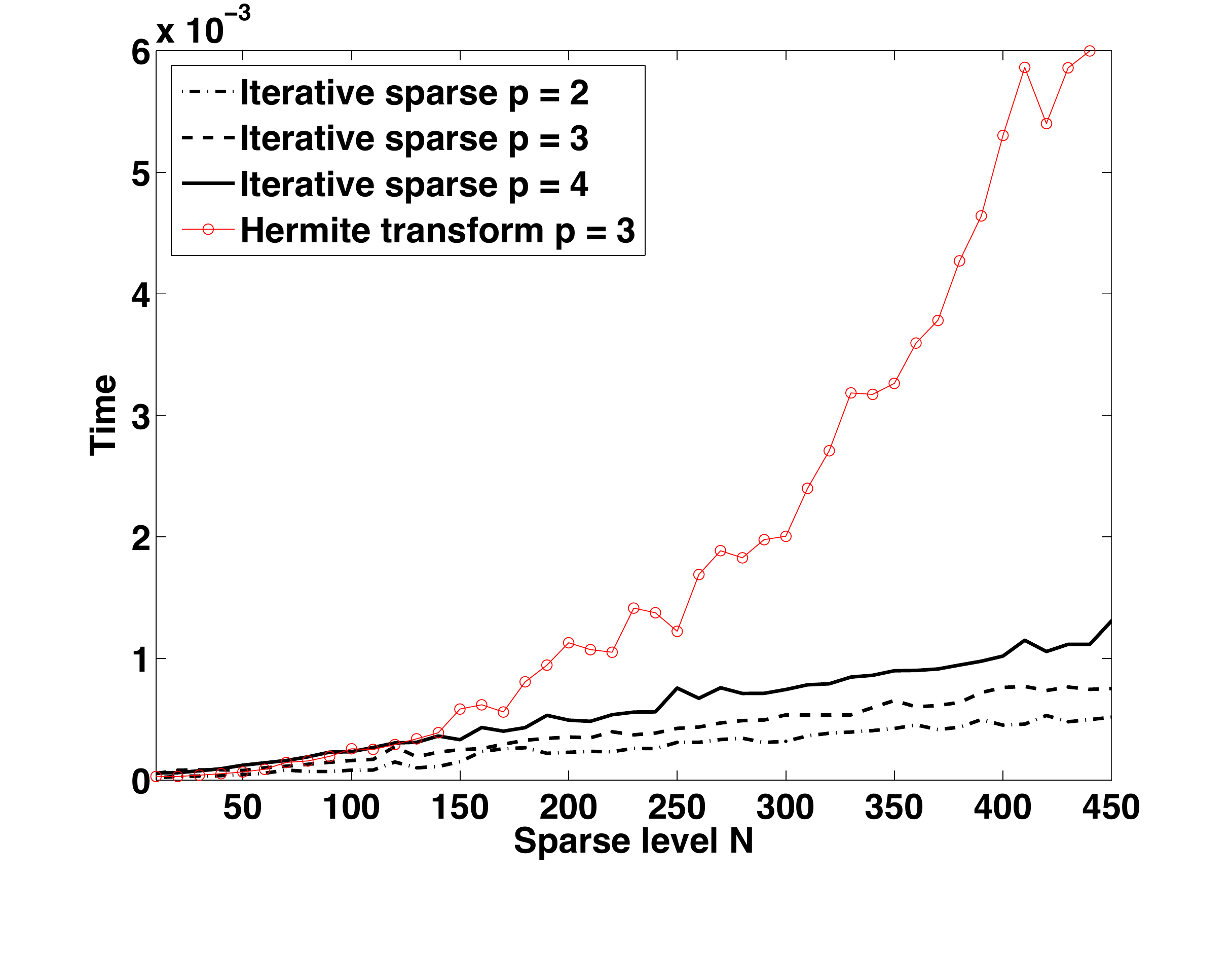}}}
\end{center}
\caption{Time versus sparse level $N$ (iterative algorithm)}
\label{fig3}
\end{figure}

In figure \ref{fig3}, we give the same time computation but with the iterative algorithm. We observe a significant speed up in the algorithm in comparison with the previous algorithm. 

\begin{figure}[ht]
\begin{center}
\rotatebox{0}{\resizebox{!}{0.4\linewidth}{%
   \includegraphics{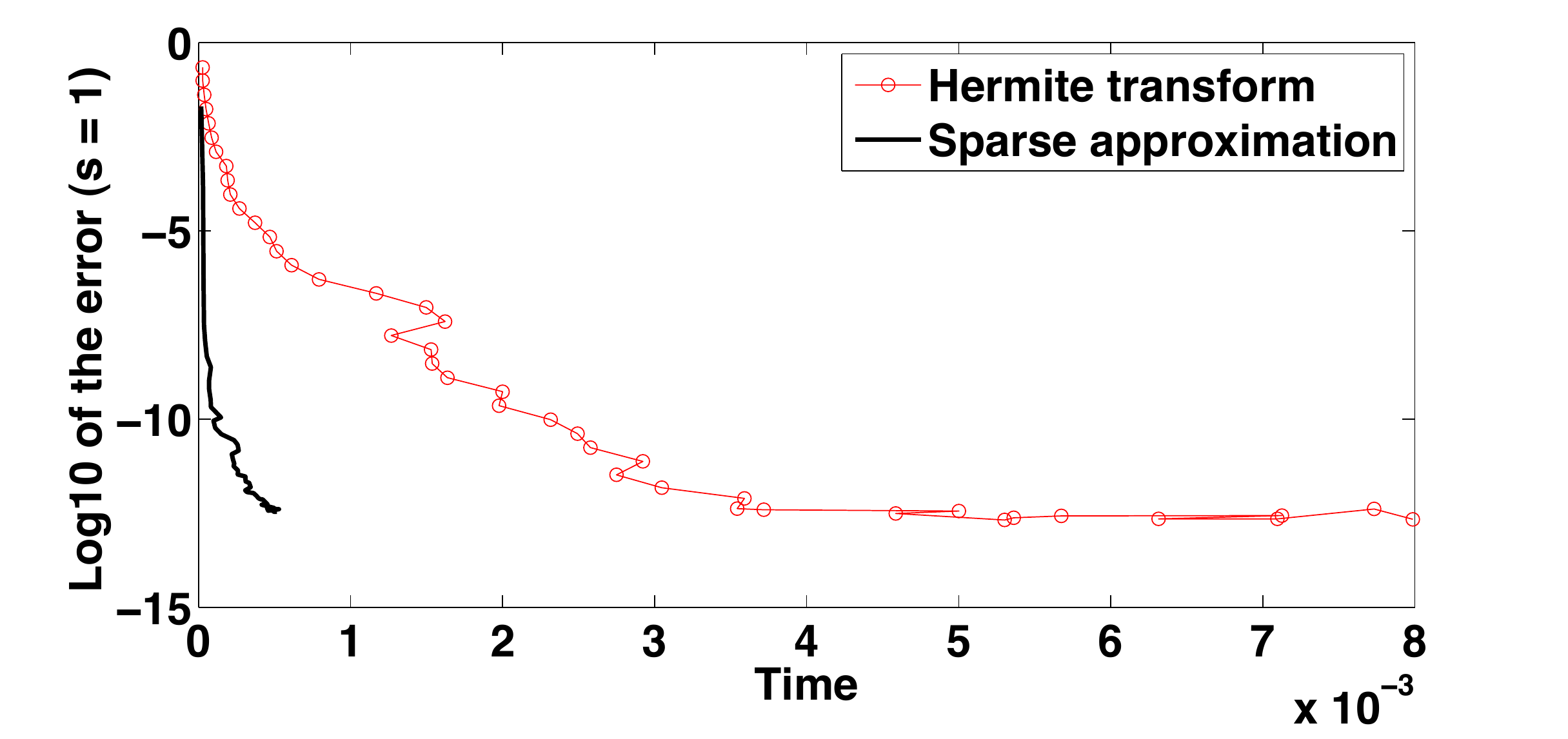}}}
\end{center}
\caption{$\sigma = 10$, $p = 2$}
\label{fig4}
\end{figure}

\begin{figure}[ht]
\begin{center}
\rotatebox{0}{\resizebox{!}{0.4\linewidth}{%
   \includegraphics{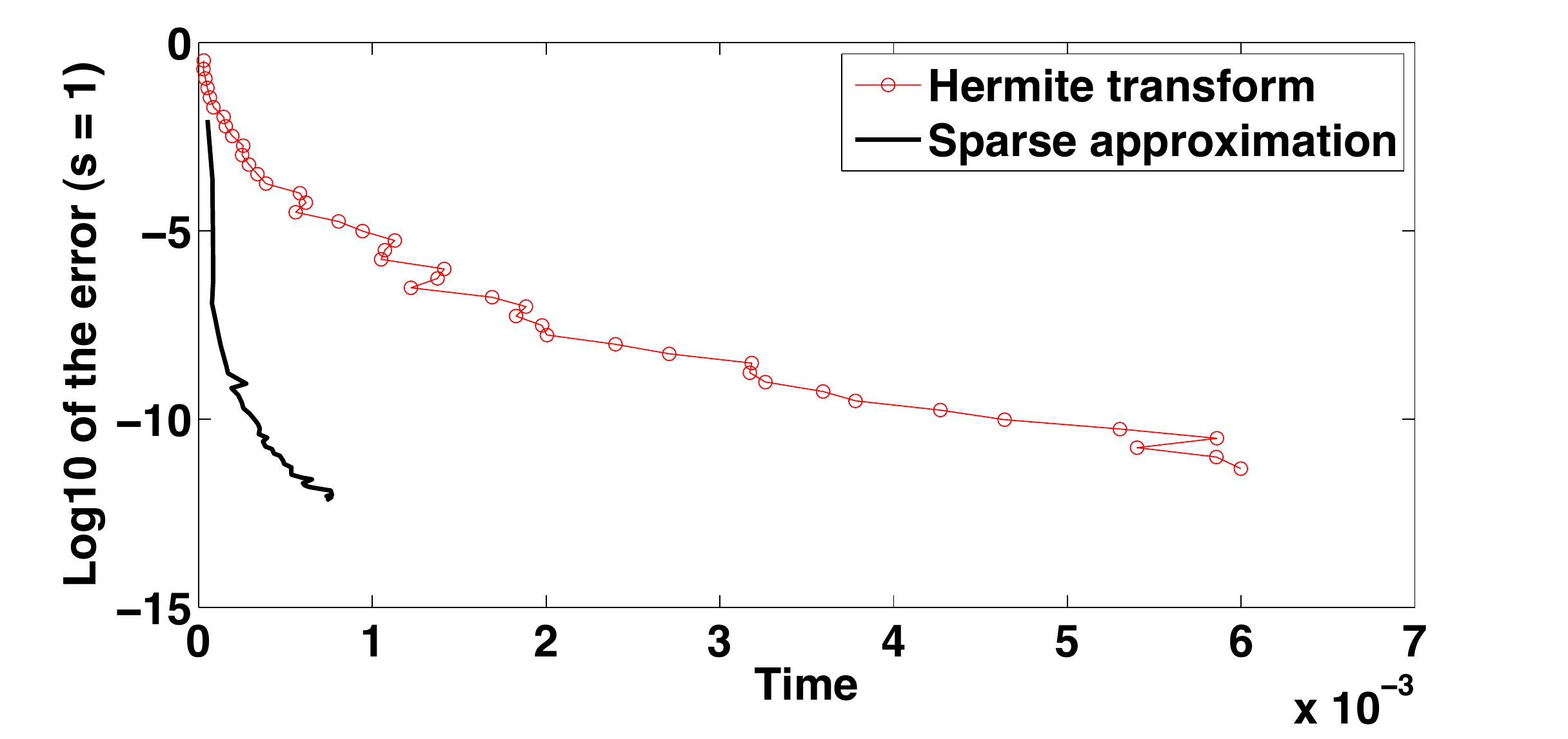}}}
\end{center}
\caption{$\sigma = 10$, $p = 3$}
\label{fig5}
\end{figure}

\begin{figure}[ht]
\begin{center}
\rotatebox{0}{\resizebox{!}{0.4\linewidth}{%
   \includegraphics{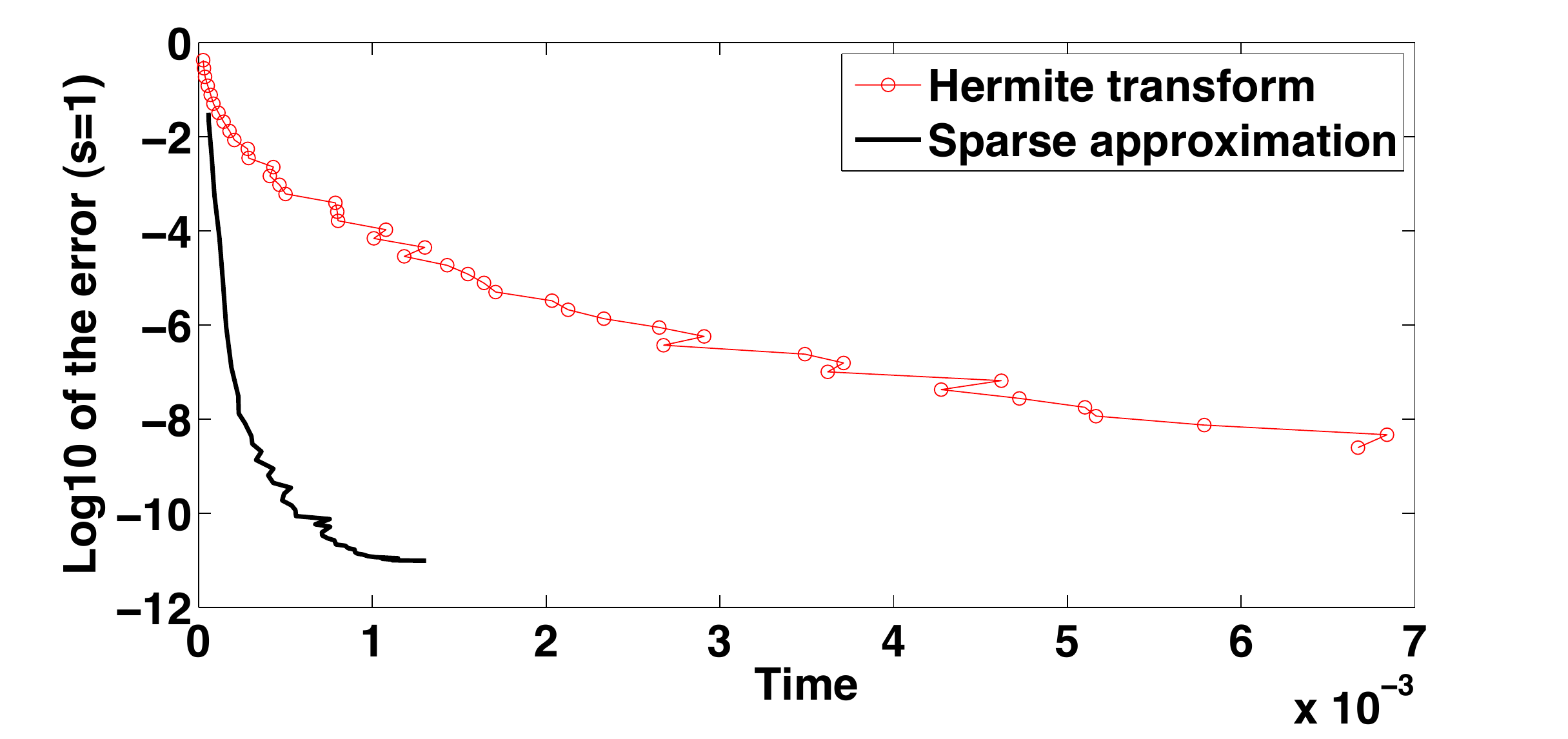}}}
\end{center}
\caption{$\sigma = 10$, $p = 4$}
\label{fig6}
\end{figure}

In the last figures \ref{fig4}, \ref{fig5} and \ref{fig6},  we fix $\sigma = 10$ for the coefficients $u_n = (1 + n)^{-\sigma}$ of the function $u(x)$, and we plot the error versus the time required for the algorithm (obtained in Figure \ref{fig2}). We compare with the result obtained with the Hermite transform method. The results obtained are better for the sparse approximation. The results obtained for the iterative algorithm are similar, but less convincing because it requires much large number $N$, despite a better cost for a single iteration. 


\section*{Appendix: Proof of Lemma \ref{lem:1}}

 To prove this Lemma, we will use the following result, which 
can be found in \cite{Greb07} for $\theta = 0$ and \cite{GIP} for $\theta = 1/2$. 

\begin{lemma}
Assume that $\jb = (j_1,\ldots,j_p) \in \Zc^p$  and for $\ell \in \Zc$, let 
\begin{equation}
\label{eq:Atheta}
A_\theta(\ell,\jb) = \frac{\mu_2(\ell,\jb)^\theta \mu_3(\ell,\jb)^{1 - \theta}}{ \mu_2(\ell,\jb)^\theta \mu_3(\ell,\jb)^{1 - \theta} + \mu_1(\ell,\jb) - \mu_2(\ell,\jb)}. 
\end{equation}
Then we have 
\begin{equation}
\label{eq:aj}
\forall\, \ell\in \Zc, \quad \Norm{\ell}{} A_{\theta}(\ell,\jb) \leq 2 \mu_1(\jb). 
\end{equation}
\end{lemma}
\begin{Proof}If $\Norm{\ell}{} \leq \mu_1(\jb)$, the equation \eqref{eq:aj} is obvious using the relation $A_{\theta}(\ell,\jb) \leq 1$. 

In the case  $\Norm{\ell}{}\geq \mu_1(\jb)$, then we have $\mu_2(\ell,\jb)^\theta \mu_3(\ell,\jb)^{1 - \theta} = \mu_1(\jb)^\theta \mu_2(\jb)^{1 - \theta}$ and 
\begin{multline*}
\Norm{\ell}{} A_{\theta}(\ell,\jb) = \frac{\Norm{\ell}{} \mu_1(\jb)^\theta \mu_2(\jb)^{1 - \theta}}{\mu_1(\jb)^\theta \mu_2(\jb)^{1 - \theta} + \Norm{\ell}{} - \mu_1(\jb)}\\
=  \left(\frac{\Norm{\ell}{} - \mu_1(\jb)}{\mu_1(\jb)^\theta \mu_2(\jb)^{1 - \theta} + \Norm{\ell}{} - \mu_1(\jb)} \right)\mu_1(\jb)^\theta \mu_2(\jb)^{1 - \theta} + \mu_1(\jb) A_{\theta}(\ell,\jb). 
\end{multline*}
We conclude using the fact that $A_{\theta}(\ell,\jb) \leq 1$, and $1 \leq \mu_1(\jb)^\theta \mu_2(\jb)^{1 - \theta} \leq \mu_1(\jb)$. 

\end{Proof}

\begin{Proofof}{Lemma \ref{lem:1}}
Let $\jb = (j_1,\ldots,j_p) \in \Zc^p$ be fixed. We distinguish three cases in the sum in $\ell \in \Zc$ appearing in \eqref{eq:lem1}. 

{\bf (i)} {$\boldsymbol{\Norm{\ell}{} > \mu_1(\jb)}$}. In this case, we have $\mu_1(\ell,\jb) = \Norm{\ell}{}$, $\mu_2(\ell,\jb) = \mu_1(\jb)$ and $\mu_{3}(\ell,\jb) = \mu_2(\jb)$. 
Hence we can write using \eqref{eq:bound} with $R = r + \kappa$, and the previous Lemma
\begin{eqnarray*}
\sum_{\Norm{\ell}{} > \mu_1(\jb)} \Norm{\ell}{}^r |a_{\ell;\jb}| &\leq& 
c \mu_2(\jb)^\nu \sum_{\Norm{\ell}{} > \mu_1(\jb)} \Norm{\ell}{}^r A_\theta(\ell,\jb)^{r + \kappa} \\
&\leq& c 2^r \mu_1(\jb)^{r}\mu_2(\jb)^\nu \sum_{\Norm{\ell}{}> \mu_1(\jb)}  A_\theta(\ell,\jb)^\kappa. 
\end{eqnarray*}
As in this case we have 
$$
A_\theta(\ell,\jb) = \frac{\mu_1(\jb)^\theta \mu_2(\jb)^{1 - \theta}}{ \mu_1(\jb)^\theta \mu_2(\jb)^{1 - \theta} + \Norm{\ell}{} - \mu_1(\jb)} \leq \mu_1(\jb)^\theta \mu_2(\jb)^{1 - \theta}\frac{1}{1 + \Norm{\ell}{} - \mu_1(\jb)}, 
$$
we obtain 
\begin{eqnarray*}
\sum_{\Norm{\ell}{}> \mu_1(\jb)} \Norm{\ell}{}^r |a_{\ell;\jb}| 
&\leq& c 2^r \mu_1(\jb)^{r + \theta \kappa}\mu_2(\jb)^{(1 - \theta) \kappa + \nu} \sum_{\Norm{\ell}{} > \mu_1(\jb)} \Big(\frac{1}{1 + \Norm{\ell}{} - \mu_1(\jb)}\Big)^\kappa\\
&\leq & C\mu_1(\jb)^{r + \theta \kappa}\mu_2(\jb)^{(1 - \theta) \kappa + \nu} , 
\end{eqnarray*}
where the constant $C$ depends only on $r$ and $\kappa$ but not on $\jb$. Here we use the fact that the sum in the right-hand side is convergent and independent of $\mu_1(\jb)$, owing to the condition $\kappa > d$.

{\bf (ii)} {$\mu_1(\jb) \geq \Norm{\ell}{} > \mu_2(\jb)$}. In this case we have $\mu_1(\ell,\jb) = \mu_1(\jb)$, $\mu_2(\ell,\jb) = \Norm{\ell}{}$ and $\mu_3(\ell,\jb) = \mu_2(\jb)$. 
Hence we have 
\begin{eqnarray*}
A_\theta(\ell,\jb) &=& \frac{\Norm{\ell}{}^\theta \mu_2(\jb)^{1 - \theta}}{\Norm{\ell}{}^\theta \mu_2(\jb)^{1 - \theta} + \mu_1(\jb) - \Norm{\ell}{}} \\
&\leq & \mu_1(\jb)^{\theta} \mu_2(\jb)^{1 - \theta} \frac{1}{1 + \mu_1(\jb) - \Norm{\ell}{}}. 
\end{eqnarray*}

Using again the previous Lemma, we obtain 
\begin{multline*}
\sum_{\mu_1(\jb) \geq \Norm{\ell}{} > \mu_2(\jb)} \Norm{\ell}{}^r |a_{\ell;\jb}| 
\leq c 2^r \mu_1(\jb)^{r}\mu_2(\jb)^\nu \sum_{\mu_1(\jb) \geq \Norm{\ell}{} > \mu_2(\jb)}  A_\theta(\ell,\jb)^\kappa\\
\leq  c 2^r \mu_1(\jb)^{r + \theta \kappa}\mu_2(\jb)^{(1 - \theta) \kappa + \nu} \sum_{\mu_1(\jb) \geq \Norm{\ell}{} > \mu_2(\jb)} \Big(  \frac{1}{1 + \mu_1(\jb) - \Norm{\ell}{}}\Big)^\kappa, 
\end{multline*}
and we conclude as in the previous case. 

{\bf (iii)} {$\mu_2(\jb) \geq \Norm{\ell}{}$}. In this last situation, we can estimate directly the term, and obtain using the fact that $A(\ell,\jb) \leq 1$, 
$$
\sum_{\Norm{\ell}{}  \leq \mu_2(\jb)} \Norm{\ell}{}^r |a_{\ell;\jb}| \leq \mu_2(\jb)^r \sum_{\Norm{\ell}{} \leq \mu_2(\jb)} |a_{\ell;\jb}|\leq 
 c \mu_2(\jb)^{\nu + r}\Big( \sum_{\Norm{\ell}{} \leq \mu_2(\jb)} 1 \Big). 
$$
Hence, using 
$\sharp \{ \ell\,\in  \Z^d\, | \, \Norm{\ell}{} \leq \mu_2(\jb)\} \leq C \mu_2(\jb)^d$, we get
$$
\sum_{\Norm{\ell}{}  \leq \mu_2(\jb)} \Norm{\ell}{}^r |a_{\ell;\jb}| \leq C \mu_2(\jb)^{r + d + \nu} \leq C \mu_1(\jb)^{r + \theta \kappa} \mu_2(\jb)^{(1 - \theta)\kappa + \nu}, 
$$
as $\kappa > d$. Gathering the previous estimate yields the result.

\end{Proofof}


\noindent {\bf Authors addresses:} 

\vskip 2ex

\noindent E. Faou, INRIA and ENS Cachan Bretagne, Avenue Robert Schumann, F-35170 Bruz, France.

\noindent {\tt Erwan.Faou@inria.fr}

\vskip 2ex 

\noindent F. Nobile, Ecole Polytechnique F\'ed\'erale de Lausanne, EPFL SB MATHICSE CSQI, MA B2 444 (B\^atiment MA), Station 8, CH-1015 Lausanne, Switzerland. 

\noindent {\tt fabio.nobile@epfl.ch}

\vskip 2ex

\noindent C. Vuillot, ENS Cachan Bretagne, Avenue Robert Schumann, F-35170 Bruz, France.

\noindent {\tt christophe.vuillot@eleves.bretagne.ens-cachan.fr}


\begin{thebibliography}{10}



\bibitem{Bam03}
D.~Bambusi, \emph{Birkhoff normal form for some nonlinear {PDE}s}, Comm. Math.
  Physics \textbf{234} (2003) 253--283.

\bibitem{Zoll}
{\rm D. Bambusi, J.-M. Delort, B. Gr\'ebert and J. Szeftel}, 
{\em Almost global existence for Hamiltonian semi-linear Klein-Gordon equations with small Cauchy data on Zoll manifolds}, Comm. Pure. Appl. Math. 60 (2007) 1665--1690.



\bibitem{BG06}
{\rm D. Bambusi and B. Gr\'ebert},
{\em Birkhoff normal form for PDE's with tame modulus}. Duke Math. J.  135  no. 3 (2006) 507–-567.

\bibitem{Griebel}
{\rm H.-J. Bungartz and M. Griebel}, {\em Sparse grids}, 
Acta Numerica (2004), pp. 1--123

\bibitem{Delort1}
{\rm J.-M. Delort and J. Szeftel}, {\em Long time existence for small data nonlinear Klein-Gordon equations on tori and spheres}.
Int. Math. Res. Not. 37 (2004) 1897--1966. 

\bibitem{Delort2}
{\rm J.-M. Delort and J. Szeftel}, {\em 
Long-time existence for semi-linear Klein-Gordon equations with small Cauchy data on Zoll manifolds }, 
Amer. J. Math. 128 (2008) 1187--1218. 


\bibitem{Greb07}
{\rm B. Gr\'ebert},  
{\em Birkhoff normal form and Hamiltonian PDEs.}
S\'eminaires et Congr\`es 15 (2007) 1--46



\bibitem{GIP}
{\rm B. Gr\'ebert, R. Imekraz and E. Paturel}, {\em 
Normal Forms for Semilinear Quantum Harmonic Oscillators}, Commun. Math. Phys. 291  (2009) 763--798.

\bibitem{Hel84}
B. Helffer, {\em Th\'eorie spectrale pour des op\'erateurs globalement elliptiques,}
Ast\'erisque, vol. 112, Soci\'et\'e Math\'ematique de France, Paris, 1984, With an English
summary.

\bibitem{Arieh}
{\rm A. Iserles}, {\em A fast and simple algorithm for the computation of Legendre coefficients}, Numer. Math. 117 (2011), 529--553. 

\bibitem{Zenger}
{\rm C. Zenger},  {\em Sparse grids, in Parallel Algorithms for Partial Differential Equations},  (W. Hackbusch, ed.), Vol. 31 of Notes on Numerical Fluid Mechanics, Vieweg, Braunschweig/Wiesbaden (1991). 

\end{thebibliography}
\end{document}